\documentclass[12pt]{article}
\usepackage{amssymb}
\usepackage{amsfonts}
\usepackage{amsmath}
\usepackage{amsthm}
\usepackage{algcompatible}
\usepackage{algorithm}
\usepackage{multirow}
\usepackage{color}

\newtheorem{theorem}{Theorem}

\newtheorem{proposition}{Proposition}

\theoremstyle{definition}

\newtheorem{example}{Example}

\theoremstyle{lemma}

\begin{document}

{\centerline{{\Large {\bf Total positivity and least squares problems}}}

\bigskip

{\centerline{{\Large {\bf in the Lagrange basis}}}

\bigskip
\bigskip
{\centerline{{\bf Ana Marco$^1$, Jos\'e-Javier Mart{\'\i}nez$^1$, Raquel Via\~na$^1$}}

\bigskip

\noindent $^1$ {\small{\it Universidad de Alcal\'a, Departamento de F{\'\i}sica y Matem\'aticas, Alcal\'a de Henares, Madrid 28871, Spain}}

\bigskip
\noindent {\small{{\bf Corresponding author:}  A. Marco. {\it E-mail:} ana.marco@uah.es}

\bigskip
\bigskip
\bigskip

\noindent{\bf Summary}

\bigskip

The problem of polynomial least squares fitting in the standard Lagrange basis is addressed in this work. Although the matrices involved in the corresponding overdetermined linear systems are not totally positive, rectangular totally positive Lagrange-Vandermonde matrices are used to take advantage of total positivity in the construction of accurate algorithms to solve the considered problem. In particular, a fast and accurate algorithm to compute the bidiagonal decomposition of such rectangular totally positive matrices is crucial to solve the problem. This algorithm also allows the accurate computation of the Moore-Penrose inverse and the projection matrix of the collocation matrices involved in these problems. Numerical experiments showing the good behaviour of the proposed algorithms are included.

\noindent

\bigskip


\bigskip
\noindent{\bf MSC:} 65F20; 15A23; 15B48
%%%% 65-XX: Numerical analysis
%%%
%%%% 65F-XX: Numerical linear algebra
%%%% 65F05: Direct methods for linear systems and matrix inversion
%%%% 65F15: Eigenvalues, eigenvectors
%%%% 65F20: Overdetermined systems, pseudoinverses
%%%
%%%% 65D-XX: Numerical approximation and computational geometry (primarily algorithms)
%%%% 65D05: Interpolation
%%%% 65D17: Computer aided design (modeling of curves and surfaces)
%%%
%%%% 15-XX: Linear and multilinear algebra; matrix theory
%%%% 15A-XX: Basic linear algebra
%%%% 15A23: Factorization of matrices
%%%
%%%% 15B-XX: Special matrices
%%%% 15B48: Positive matrices and their generalizations; cones of matrices

\bigskip
\bigskip
\bigskip

\section{Introduction} %%%%%%%%%%%%%%%%%%%%%%%%%%%%%%%%%%%%%%%%%%%%%%%%%%%%%%%%%%%%%%%%%%%%%%%%%%%%%%%%%%%%%%%%%%%%%%%
\label{sec:Intro}

Given the nodes $x_j\in\mathbb{R}$, for $j=1,\ldots, n+1$, the associated {\em Lagrange basis} of the space $\Pi_n(t)$ of the polynomials of degree less than or equal to $n$ is

\begin{equation}
\label{eq:basisLagrange}
\mathcal{B}_L=\left\{ \prod_{\substack{k=1\\ k \neq 1}}^{n+1}\frac{t-x_k}{x_1-x_k},~ \prod_{\substack{k=1\\ k \neq 2}}^{n+1}\frac{t-x_k}{x_2-x_k}, \ldots , %\prod_{i=1,~i \neq n+1}^{n+1}\frac{t-x_i}{x_{n+1}-x_i}\right\}.
\prod_{\substack{k=1\\ k \neq n+1}}^{n+1}\frac{t-x_k}{x_{n+1}-x_k}\right\}.
\end{equation}
This basis is widely used in polynomial interpolation \cite{BT04,CF13,CW04,T13}, as the data $(x_j,y_j)$, for $1\leq j\leq n+1$, directly provide the interpolating polynomial,
$$p(t)=\sum_{j=1}^{n+1}y_j~l_j(t),~\text{ with }~l_j(t)=\prod_{\substack{k=1\\ k \neq j}}^{n+1}\frac{t-x_k}{x_j-x_k}.$$
This representation of the interpolating polynomial is sometimes very useful, for example for computing the roots of the polynomial or for evaluating it \cite{BT04,CW04}.  

In this work a different problem involving the Lagrange basis is considered. Besides the nodes $\{x_j\}_{1\leq j \leq n+1}$, required to build the basis $\mathcal{B}_L$, a data set is provided $\{(t_i,b_i)\}_{1\leq i \leq l+1}$, with $l\geq n$. The {\it least squares fitting polynomial} of degree less than or equal to $n$ for the given data will be built, obtaining the coordinates in the Lagrange basis of such polynomial. As it is well known, solving this problem is equivalent to solve, in the least squares sense, the overdetermined linear system $Lc=b$, where $L$ is the {\it collocation matrix} corresponding to the basis $\mathcal{B}_L$ and to the points $\{t_i\}_{1\leq i\leq l+1}$ (see (\ref{eq:L})).

Since matrix $L$ is an ill-conditioned structured matrix \cite{MMVLag}, standard linear algebra algorithms to solve problems involving it provide no accurate results. Due to this fact, specific algorithms taking into account the structure of matrix $L$ must be designed. To develop them, results from the field of totally positive matrices, mainly involving the {\it bidiagonal decomposition} of such matrices, will be essential \cite{KOEV05,KOEV07}.

Classically, a matrix is called {\it totally positive} (TP) if all its minors are nonnegative and {\it strictly totally positive} (STP) if all its minors are positive \cite{Ando,GP96,PINKUS}. Although we will follow such terminology, we must note that totally positive matrices and strictly totally positive matrices are also called {\it totally nonnegative} and {\it totally positive} matrices, respectively \cite{FJ11}.

We observe that, although $L$ is not a totally positive matrix, it can be factorized as the product of a totally positive matrix and a diagonal matrix, which will allow us to employ total positivity results to solve accurately our problem. Moreover, our aim will be to solve as much stages of the problem as possible to {\it high relative accuracy}.

Let us recall that real number $x$ is computed to high relative accuracy (HRA) whenever the computed value $\widehat{x}$ satisfies
$$
\frac{|x-\widehat{x}|}{|x|}\leq K u,
$$
where $u$ is the round-off unit and $K>0$ is a constant independent of the arithmetic precision. On the other hand, an algorithm computes to HRA if it satisfies the so called  {\it no inaccurate cancellation} (NIC) condition \cite{DDHK}: The algorithm only multiplies, divides, adds (resp. subtracts) real numbers with like (resp. differing) signs, and otherwise only adds or subtracts input data.

The starting point of the work developed here is \cite{MMVLag}. In such paper Lagrange-Vandermonde square matrices are considered. They are collocation matrices corresponding to bases created by removing denominators from the standard Lagrange basis $\mathcal{B}_L$. Under certain conditions, square Lagrange-Vandermonde matrices are proven to be strictly totally positive, and this fact allows the creation of fast algorithms to solve to HRA the problems of inverse computation and linear system solving (with a data vector having an alternating sign pattern) involving the square collocation matrix $L$ associated to the standard Lagrange basis $\mathcal{B}_L$.

In this work, the problems of computing the {\it Moore-Penrose inverse} and the {\em projection matrix} of matrix $L$ in an accurate and efficient way are also analyzed. One of the main applications of the Moore-Penrose inverse is the computation of the solution of least squares problems. %Sacado de \cite{MM19}
The projection matrix is widely used in regression problems \cite{HW}, in the detection of high-leverage points or in the computation of the least squares residuals.

The rest of the paper is organized as follows. In Section 2 we prove that under certain conditions, which will be assumed to hold until Section 4 included, a rectangular Lagrange-Vandermonde matrix $A$ is strictly totally positive, and we provide its bidiagonal factorization. In Section 3, the $QR$ decomposition of $A$ is obtained, and it is used to solve the least squares problem $Lc=b$. In Section 4, the Moore-Penrose inverse and the projection matrix of $L$ are obtained. The general situation, where $A$ is not an STP matrix, is studied in Section 5. Finally, in Section 6 some numerical examples are provided.

\section{Total positivity and bidiagonal factorization of a rectangular Lagrange-Vandermonde matrix} %%%%%%%%%%%%%%%%%%%%%%%%%%%%%%%%%%%%%%%%%%%%%%%%%%%%%%%%%%%%%%%%%%%%%%%%%%%%%%%%%%%%%%%%%%%%%%%
\label{sec:TNBDLR}
In this section we extend to the rectangular case the fast and accurate algorithm to compute the bidiagonal factorization of Lagrange-Vandermonde matrices developed in the square case in Section 4 of \cite{MMVLag}.

Let us consider the basis $\mathcal{B}$ of the space $\Pi_n(t)$ of polynomials of degree less than or equal to $n$, formed by removing the denominators of the Lagrange basis (\ref{eq:basisLagrange}),

\begin{equation}
\label{eq:basis}
\mathcal{B}=\left\{ \prod_{\substack{k=1\\ k \neq 1}}^{n+1}(t-x_k), \prod_{\substack{k=1\\ k \neq 2}}^{n+1}(t-x_k), \ldots \prod_{\substack{k=1\\ k \neq n+1}}^{n+1}(t-x_k)\right\},
\end{equation}
where $x_j \in \mathbb R$, $j=1,\ldots, n+1$. The collocation matrix of this basis at the nodes $\{t_i\}_{1\leq i \leq l+1}$, with $l\geq n$, is

\begin{equation}
\label{eq:A}
A=\left(
    \begin{array}{cccc}
      \displaystyle{\prod_{\substack{k=1\\ k \neq 1}}^{n+1}(t_1-x_k)} & \displaystyle{\prod_{\substack{k=1\\ k \neq 2}}^{n+1}(t_1-x_k)} & \cdots & \displaystyle{\prod_{\substack{k=1\\ k \neq n+1}}^{n+1}(t_1-x_k)} \\
        & & & \\
       \displaystyle{\prod_{\substack{k=1\\ k \neq 1}}^{n+1}(t_2-x_k)} & \displaystyle{\prod_{\substack{k=1\\ k \neq 2}}^{n+1}(t_2-x_k)} & \cdots & \displaystyle{\prod_{\substack{k=1\\ k \neq n+1}}^{n+1}(t_2-x_k)} \\
               & & & \\
      \vdots & \vdots & \ddots & \vdots \\
       \displaystyle{\prod_{\substack{k=1\\ k \neq 1}}^{n+1}(t_{l+1}-x_k)} & \displaystyle{\prod_{\substack{k=1\\ k \neq 2}}^{n+1}(t_{l+1}-x_k)} & \cdots & \displaystyle{\prod_{\substack{k=1\\ k \neq n+1}}^{n+1}(t_{l+1}-x_k)}\\
    \end{array}
  \right).
\end{equation}

The following theorem shows that, under certain conditions, $A$ is an STP matrix and provides its bidiagonal factorization.

\begin{theorem}
\label{th:FactorA}
Let $A=(a_{i,j})_{1\le i \leq l+1;1\leq j \leq n+1}$ be a
Lagrange-Vandermonde matrix corresponding to the basis $\mathcal{B}$ and the nodes $\{t_j\}_{1\leq j \leq l+1}$, where $x_1<x_2<\ldots<x_{n+1}<t_{l+1}<\ldots<t_2<t_1$, with $l\geq n$. Then $A$ is $STP$ and it admits a factorization in
the form
\begin{equation}
\label{eq:FactorA}
A=F_{l}F_{l-1} \cdots F_1 D G_1 \cdots G_{n-1}G_{n} ,
\end{equation}
where $F_i$ ($1\leq i \leq l$) are $(l+1)\times (l+1)$  bidiagonal matrices of the form
\begin{equation}
\label{eq:Fi}
F_i=\left(
\begin{smallmatrix}
1 & & & & & & & \\
0 & 1 & & & & & & \\
& \ddots & \ddots & & & & & \\
& & 0 & 1 & & & & \\
& & & m_{i+1,1} & 1 & & & \\
& & & & m_{i+2,2} & 1 & & \\
& & & & & \ddots & \ddots & \\
& & & & & & m_{l+1,l+1-i} & 1
 \end{smallmatrix}
 \right),
 \end{equation}
$G_i^T$ ($1\leq i \leq n$) are $(n+1)\times (n+1)$  bidiagonal matrices of the form

\begin{equation}
\label{eq:GiT}
G_i^T=\left(
\begin{smallmatrix}
1 & & & & & & & \\
0 & 1 & & & & & & \\
& \ddots & \ddots & & & & & \\
& & 0 & 1 & & & & \\
& & & \widetilde m_{i+1,1} & 1 & & & \\
& & & & \widetilde m_{i+2,2} & 1 & & \\
& & & & & \ddots & \ddots & \\
& & & & & & \widetilde m_{n+1,n+1-i} & 1
 \end{smallmatrix}
 \right),
\end{equation}
and $D$ is the $(l+1)\times (n+1)$ diagonal matrix
\begin{equation}
\label{eq:D}
D =\textnormal{diag}\{p_{1,1},p_{2,2},\ldots,p_{n+1,n+1}\}.
\end{equation}

Besides, the elements $m_{i,j}$ and $\widetilde m_{i,j}$ of these matrices are the multipliers of the Neville elimination of $A$ and $A^T$ respectively, and the elements $p_{i,i}$ are the diagonal pivots of the Neville elimination of $A$. The expressions of all these elements are
\begin{equation}
\label{eq:mij}
m_{i,j}=(t_{i-j}-x_j)\frac{\displaystyle{\prod_{k=i-j+1}^{i-1}(t_i-t_k) \prod_{k=j+1}^{n+1}(t_i-t_k)}}{\displaystyle{\prod_{k=i-j}^{i-2}(t_{i-1}-t_k)\prod_{k=j}^{n+1}(x_k-t_{i-1})}}, %\eqno(4)
\end{equation}
where $j=1,\ldots,n+1$; $i=j+1, \dots, l+1$,

\begin{equation}
\label{eq:mijT}
\widetilde m_{i,j}=(t_j-x_{i-j}) \frac{\displaystyle{\prod_{k=i-j+1}^{i-1}(x_i-x_k) \prod_{k=1}^{j-1}(t_k-x_{i-1})}}{\displaystyle{\prod_{k=i-j}^{i-2}(x_{i-1}-x_k) \prod_{k=1}^{j}(t_k-x_i)}}, %\eqno(5)
\end{equation}
where $j=1,\dots,n$; $i=j+1,\ldots,n+1$, and

\begin{equation}
\label{eq:pii}
p_{i,i}=\frac{\displaystyle{\prod_{k=1}^{i-1}(t_i-t_k) \prod_{k=1}^{i-1}(x_i-x_k) \prod_{k=i+1}^{n+1}(t_i-x_k)}}{\displaystyle{\prod_{k=1}^{i-1}(t_k-x_i)}}, %\eqno(6)
\end{equation}
where $i=1,\ldots,n+1$.

\end{theorem}

\begin{proof}
In \cite{MMVLag} the strict total positivity of matrix $A$ has been proved for the $(n+1)\times (n+1)$ square case by demonstrating that each minor of A with initial consecutive rows and each minor of $A$ with initial consecutive columns is positive. The extension to the rectangular case is straightforward, and therefore we can assure that rectangular matrix $A$ is STP.

According to \cite{GP96}, matrix $A$ admits the factorization (\ref{eq:FactorA}), with the multiplying matrices given by (\ref{eq:Fi}), (\ref{eq:GiT}) and (\ref{eq:D}). A direct extension to the rectangular case of the formulae obtained in \cite{MMVLag} for the multipliers of the Neville elimination of $A$ and $A^T$ respectively, and for the diagonal pivots of the Neville elimination of $A$ provides expressions (\ref{eq:mij}), (\ref{eq:mijT}), and (\ref{eq:pii}).
\end{proof}

Following the terminology of \cite{KOEV05,KOEV07} we call $\mathcal{BD}(A)$ to the matrix of size $(l+1)\times(n+1)$ containing the bidiagonal decomposition of the Lagrange-Vandermonde matrix $A$ as follows

$$\mathcal{BD}(A)=\left(
\begin{array}{cccc}
p_{11}    & \tilde{m}_{21} & \cdots & \tilde{m}_{n+1,1} \\
m_{21}    & p_{22}         & \cdots & \tilde{m}_{n+1,2} \\
\vdots    & \vdots         & \ddots & \vdots \\
m_{n+1,1} & m_{n+1,2}      & \cdots & p_{n+1,n+1} \\
\vdots    & \vdots         & \ddots & \vdots \\
m_{l+1,1} & m_{l+1,1}      & \cdots & m_{l+1,n+1}
\end{array}
\right),$$
where $m_{ij}$ are the multipliers of the Neville elimination of $A$, $\tilde{m}_{ij}$ are the multipliers of the Neville elimination of $A^T$ and $p_{ii}$ are the diagonal pivots of the Neville elimination of $A$. Their explicit expressions (\ref{eq:mij}), (\ref{eq:mijT}) and (\ref{eq:pii}) respectively provide the bidiagonal factorization of $A$ to HRA, which cannot be achieved by performing directly the Neville elimination of $A$ \cite{KOEV05}.

The extension to the rectangular case of the algorithm {\tt TNBDL} presented in \cite{MMVLag} provides an algorithm, which will be called {\tt TNBDLR}\footnote{This and the rest of the algorithms presented in this paper are available under request.}, to compute the matrix $\mathcal{BD}(A)$, where $A$ is an STP rectangular Lagrange-Vandermonde matrix, for $\{x_j\}_{1\leq j \leq n+1}$ and $\{t_i\}_{1\leq i \leq l+1}$. The algorithm is fast because it requires $O(ln)$ arithmetic operations, and it has high relative accuracy because it does not perform any substractive cancellation operation, hence satisfying the NIC condition. In addition, the algorithm does not require the explicit construction of matrix $A$, with the corresponding storage space saving.

The $\mathcal{BD}(A)$, computed to HRA as previously explained, will be the starting point of the algorithms proposed in the following sections.

\section{Least squares problem} %%%%%%%%%%%%%%%%%%%%%%%%%%%%%%%%%%%%%%%%%%%%%%%%%%%%%%%%%%%%%%%%%%%%%%%%%%%%%%%%%%%%%%%%%%%%%%%
\label{sec:LS}

In this section we address the problem of, given data $\{(t_i,b_i)\}_{1\leq i\leq l+1}$, finding the $n$-degree  (with $l\geq n$) polynomial

\begin{equation}
\label{eq:polLagrange}
p(t)=\sum_{j=1}^{n+1}c_j~l_j(t),~\text{ with }~l_j(t)=\prod_{\substack{k=1\\ k \neq j}}^{n+1}\frac{t-x_k}{x_j-x_k},
\end{equation}
expressed in the Lagrange basis $\mathcal{B}_L$ (\ref{eq:basisLagrange})
 which best fits in the least squares sense the data, that is, which minimizes $$\sum_{i=1}^{l+1}~\left|b_i-p(t_i)\right|^2.$$

The unknowns of our problem are the coefficients $c_j$ of the polynomial (\ref{eq:polLagrange}). Thus, solving this problem is equivalent to solve in the least squares sense the overdetermined linear system

$$Lc=b,$$
where $L$ is the collocation matrix corresponding to the base $\mathcal{B}_L$ and to the nodes $\{t_i\}_{1\leq i \leq l+1}$,
{\small\begin{equation}L=\left(\begin{array}{cccc}
\displaystyle\prod_{\substack{k=1\\ k \neq 1}}^{n+1}\frac{t_1-x_k}{x_1-x_k} & \displaystyle\prod_{\substack{k=1\\ k \neq 2}}^{n+1}\frac{t_1-x_k}{x_2-x_k} & \ldots & \displaystyle\prod_{\substack{k=1\\ k \neq n+1}}^{n+1}\frac{t_1-x_k}{x_{n+1}-x_k}\\[8mm]
\displaystyle\prod_{\substack{k=1\\ k \neq 1}}^{n+1}\frac{t_2-x_k}{x_1-x_k} & \displaystyle\prod_{\substack{k=1\\ k \neq 2}}^{n+1}\frac{t_2-x_k}{x_2-x_k} & \ldots & \displaystyle\prod_{\substack{k=1\\ k \neq n+1}}^{n+1}\frac{t_2-x_k}{x_{n+1}-x_k}\\[2mm]
\vdots & \vdots & \ddots & \vdots \\[2mm]
\displaystyle\prod_{\substack{k=1\\ k \neq 1}}^{n+1}\frac{t_{l+1}-x_k}{x_1-x_k} & \displaystyle\prod_{\substack{k=1\\ k \neq 2}}^{n+1}\frac{t_{l+1}-x_k}{x_2-x_k} & \ldots & \displaystyle\prod_{\substack{k=1\\ k \neq n+1}}^{n+1}\frac{t_{l+1}-x_k}{x_{n+1}-x_k}
\end{array}\right).\label{eq:L}\end{equation}}

Let us notice that matrix $L$ can be factorized as

\begin{equation}L=A\overline{D},\label{eq:LAD}\end{equation}
where $\overline D$ is the diagonal matrix
\begin{equation}
\label{eq:OverlineD}\overline{D}=\left(\begin{array}{cccc}
\displaystyle\prod_{\substack{k=1\\ k \neq 1}}^{n+1}\frac{1}{x_1-x_k} & 0 & \cdots & 0 \\
0 & \displaystyle\prod_{\substack{k=1\\ k \neq 2}}^{n+1}\frac{1}{x_2-x_k} & \cdots & 0 \\
\vdots & \vdots & \ddots & \vdots \\[2mm]
0 & 0 & \cdots & \displaystyle\prod_{\substack{k=1\\ k \neq n+1}}^{n+1}\frac{1}{x_{n+1}-x_k}
\end{array}\right),\end{equation}
and $A$ is the STP matrix in (\ref{eq:A}).

In the next proposition we show how the solutions in the least squares sense of the linear systems $Lc=b$ and $Az=b$ are related.

\begin{proposition}
\label{prop:SolLS}
Given the nodes $\{x_j\}_{1\leq j\leq n+1}$ and the data $\{(t_i,b_i)\}_{1\leq i\leq l+1}$, %such that the order given in  (\ref{eq:NodesOrder}) is satisfied,
let $\overline z$ be the unique least-squares solution of the linear system $Az=b$, where $A$ is the matrix in (\ref{eq:A}). Then, the unique solution in the least squares sense of the linear system $Lc=b$, where $L$ is the matrix in (\ref{eq:L}), is $\overline c=\overline D ^{-1}\overline z$, with $\overline D$ being the diagonal matrix of (\ref{eq:OverlineD}).
\end{proposition}

\begin{proof}
Since matrix $L$ has full rank, the least squares solution of the overdetermined linear system $Lc=b$ is given by the unique solution of the normal equations, $L^TL\overline c=L^Tb$. Taking into account that $L=A\overline D$, these equations can be written as $\overline D^TA^TA\overline D\overline c=\overline D^TA^Tb$. As $\overline D$ is a square non-singular matrix, the above equations are equivalent to $A^TA\overline D\overline c=A^T b$. Therefore, calling $\overline z=\overline D\overline c$ it is satisfied that $\overline z$ is the least-squares solution of the system $Az=b$.
\end{proof}

According to the previous proposition, solving the least squares system $Lc=b$ can be decomposed into two steps: first, solve the least squares system $Az=b$; second, compute $\overline c=\overline D^{-1}\overline z$, with $\overline D$ given in (\ref{eq:OverlineD}). Algorithm \ref{alg:LSLagTypeMMV}, which will be explained in detail below, performs the first step and Algorithm \ref{alg:LSLagMMV} provides the least squares solution of the system $Lc=b$.

From now on, we consider the following ordering
\begin{equation}x_1<x_2<\ldots<x_{n+1}<t_{l+1}<\ldots<t_2<t_1,\label{eq:NodesOrder}\end{equation}
the more general setting being considered in Section \ref{sec:GeneralCase}.

The matrix $L$ is not TP but, provided (\ref{eq:NodesOrder}) is satisfied, the corresponding Lagrange-Vandermonde matrix $A$ is STP (see Theorem \ref{th:FactorA}). The strict total positivity of $A$ will be the key to develop a fast and accurate algorithm to solve the least squares problem at hand.

The fact that $A$ is a full rank STP matrix makes the method based on the QR decomposition due to Golub \cite{GOL} adequate to solve the first step \cite{BJO}. A useful presentation of important results on the numerical solution of least squares problems can be found in Chapter 5 of \cite{I09}.

For the sake of completeness, we include the following result (see Section 1.3.1 in \cite{BJO}).

\medskip
\begin{theorem}
\label{tm:QR}
Let $Ax=b$ be a linear system, where $A \in \mathbb{R}^{(l+1) \times (n+1)}$, $l \geq n$, $x \in \mathbb{R}^{n+1}$ and $b \in \mathbb{R}^{l+1}$. Assume that $rank(A)=n+1$, and let the QR decomposition of $A$  be given by
\begin{equation}
\label{eq:QR}
A=Q \left[\begin{array}{c} R\\ 0 \end{array}\right],
\end{equation}
where $Q \in \mathbb{R}^{(l+1) \times (l+1)}$ is an orthogonal matrix and $R \in \mathbb{R}^{(n+1) \times (n+1)}$ is an upper triangular matrix with positive diagonal entries. Then the solution of the least squares problem $$min_x \parallel \hspace*{-0.1cm}b-Ax \hspace*{-0.1cm}\parallel_2$$ is obtained from
$$
\left[\begin{array}{c} d_1 \\ d_2\end{array}\right] = Q^T b, \quad Rx=d_1, \quad r=Q \left[\begin{array}{c}0 \\ d_2 \end{array}\right],
$$
where $d_1 \in \mathbb{R}^{n+1}$, $d_2\in \mathbb{R}^{l-n}$ and $r=b-Ax$. In particular $\parallel \hspace*{-0.1cm}r\hspace*{-0.1cm} \parallel_2 = \parallel \hspace*{-0.1cm}d_2\hspace*{-0.1cm} \parallel_2$.
\end{theorem}

Taking into account this result, the least squares solution of the system $Az=b$ is computed in  Algorithm \ref{alg:LSLagTypeMMV}.

\begin{algorithm}[!h]
\caption{}
\begin{algorithmic}[1]
\REQUIRE Vectors $x\in \mathbb{R}^{n+1}$, $t\in \mathbb{R}^{l+1}$ and $b\in \mathbb{R}^{l+1}$ containing the nodes $\{x_j\}_{1\leq j\leq n+1}$ and

the data $\{(t_i,b_i)\}_{1\leq i\leq l+1}$ satisfying (\ref{eq:NodesOrder}).
\ENSURE The solution vector $\overline z\in \mathbb{R}^{n+1}$, in the least squares sense, of the linear system

$Az=b$.

\medskip
\STATE function $\overline z$={\tt LSSolveA}($x$,$t$,$b$)

\medskip
\STATE $B$={\tt TNBDLR}($x$,$t$);

\medskip
\STATE $[Q,C]={\tt TNQR}(B)$;
\STATE $C_1=C(1:n+1,:)$;
\STATE $Q_1=Q(:,1:n+1)$;
\STATE $d_1=Q_1'*b'$;
\STATE $\overline z={\tt TNSolve}(C_1,d_1)$;
\end{algorithmic}\label{alg:LSLagTypeMMV}
\end{algorithm}

The algorithm {\tt TNQR} has been developed by P. Koev, and given the bidiagonal factorization of $A$, it computes the matrix $Q$ and the bidiagonal factorization of the matrix $R$. Let us point out here that if $A$ is STP, then $R$ is TP. {\tt TNQR} has a computational cost of $O(l^2n)$ arithmetic operations if the matrix $Q$ is required \cite{KOEV07}, and its implementation in \textsc{Matlab} can be obtained from \cite{KOEV}.

Line $6$ of Algorithm \ref{alg:LSLagTypeMMV} is carried out by using the standard matrix multiplication command of \textsc{Matlab}. As for line $7$, it solves the system $C_1z=d_1$ by means of the algorithm {\tt TNSolve} of Koev \cite{KOEV}. It has a computational cost of $O(n^2)$ arithmetic operations, and its HRA is only guaranteed when $b$ has an alternating sign pattern.

Taking into account the previous explanation, the total cost of Algorithm \ref{alg:LSLagTypeMMV} is $O(l^2n)$ arithmetic operations.

\begin{algorithm}[!h]
\caption{}
\begin{algorithmic}[1]
\REQUIRE Vectors $x\in \mathbb{R}^{n+1}$, $t\in \mathbb{R}^{l+1}$ and $b\in \mathbb{R}^{l+1}$ containing the nodes $\{x_j\}_{1\leq j\leq n+1}$ and

the data $\{(t_i,b_i)\}_{1\leq i\leq l+1}$ satisfying (\ref{eq:NodesOrder}).
\ENSURE The solution vector $\overline c\in \mathbb{R}^{n+1}$, in the least squares sense, of the linear system

$Lc=b$.

\medskip
\STATE function ${\overline c}$= {\tt LSSolveL}($x$,$t$,$b$)

\medskip
\STATE $z$= {\tt LSSolveA}($x$,$t$,$b$)

\medskip
\STATE {\bf for} $j=1:n+1$
%\STATE \quad $den_j=\displaystyle\prod_{\substack{k=1\\ k \neq j}}^{n+1}(x_j-x_k)$;
\STATE \quad $den(j)={\tt prod}(x(j)-x(1:j-1))*{\tt prod}(x(j)-x(j+1:n+1))$;
\STATE \quad $\overline c(j)=den(j)* z(j)$;
\STATE {\bf end}

\medskip
\STATE $\overline c$;
\end{algorithmic}\label{alg:LSLagMMV}
\end{algorithm}

In Algorithm \ref{alg:LSLagMMV}, the product $\overline D ^{-1}\overline z$ is performed in lines 3-6. The computational cost of the algorithm is led by the cost of computing the QR decomposition of $A$ in Algorithm \ref{alg:LSLagTypeMMV}, and therefore, it is of $O(l^2n)$ arithmetic operations.

\section{Moore-Penrose inverse and projection matrix}
\label{sec:MPInverseProjMatrix}

Our aim in this section is to compute in an accurate and efficient way the Moore-Penrose inverse and the projection matrix of the collocation matrix $L$ in (\ref{eq:L}), corresponding to the Lagrange basis $\mathcal{B}_L$ of (\ref{eq:basisLagrange}) and the nodes $\{t_i\}_{1\leq i\leq l+1}$, where the ordering given in (\ref{eq:NodesOrder}) is satisfied.

Let us recall that the {\it Moore-Penrose inverse} of a matrix $N\in \mathbb{R}^{(l+1)\times (n+1)}$, usually denoted by $N^{\dagger}$, is the unique matrix $G\in \mathbb{R}^{(n+1)\times (l+1)}$ satisfying the four {\it Penrose conditions}:
\begin{enumerate}
\item $NGN=N$, \item $GNG=G$, \item $(NG)^T=NG$, \item $(GN)^T=GN$.
\end{enumerate}

It can be easily checked that, if $l\geq n$ and rank($N$)=$n$, it is satisfied that $$N^{\dagger}=(N^TN)^{-1}N^T.$$

In this situation, Theorem 3.1 of \cite{MM19} is satisfied, and therefore $N^{\dagger}$ can be computed from the $QR$ decomposition (\ref{eq:QR}) of $N$ as $$N^{\dagger}=R^{-1}Q_1^T,$$ where $Q_1$ is the $(l+1)\times (n+1)$ matrix with the first $n+1$ columns of $Q$.

\begin{proposition}
Given the nodes $\{x_j\}_{1\leq j\leq n+1}$ and $\{t_i\}_{1\leq i\leq l+1}$ satisfying  (\ref{eq:NodesOrder}),
 the Moore-Penrose inverse $L^{\dagger}$ of matrix $L$ in (\ref{eq:L}) is $$L^{\dagger}=\overline D^{-1}A^{\dagger},$$ where $\overline D$ is the diagonal matrix in (\ref{eq:OverlineD}) and $A$ is the Lagrange-Vandermonde matrix of (\ref{eq:A}). \label{prop:MPInverse}
\end{proposition}

\begin{proof}
Since $L$ has full rank, its Moore-Penrose inverse is $L^{\dagger}=(L^TL)^{-1}L^T$. From the factorization $L=A\overline D$ given in (\ref{eq:LAD}), taking into account that $A$ is a full rank matrix and $\overline D$ is a nonsingular matrix, it follows that
$$L^{\dagger}=%((A\overline D)^T(A\overline D))^{-1}(A\overline D)^T=(\overline D^TA^TA\overline D)^{-1}\overline D^TA^T=
\overline D^{-1}\Big(A^TA\Big)^{-1}\Big(\overline D^T\Big)^{^{-1}}\overline D^TA^T=\overline D^{-1}\Big(A^TA\Big)^{-1}A^T=\overline D^{-1}A^{\dagger}.$$
\end{proof}

Using this result and the fact that matrix $A$ is STP (see Theorem \ref{th:FactorA}), the Moore-Penrose inverse $L^{\dagger}$ can be accurately computed using Algorithm \ref{alg:PseudoInverseL}, which calls to Algorithm \ref{alg:PseudoInverseA} to perform the Moore-Penrose inverse $A^{\dagger}$ of $A$. 

\begin{algorithm}[!h]
\caption{}
\begin{algorithmic}[1]
\REQUIRE Vectors $x\in \mathbb{R}^{n+1}$ and $t\in \mathbb{R}^{l+1}$ containing the nodes $\{x_j\}_{1\leq j\leq n+1}$ and $\{t_i\}_{1\leq i\leq l+1}$

satisfying (\ref{eq:NodesOrder}).
\ENSURE The Moore-Penrose inverse $A^{\dagger}$ of matrix $A$ in (\ref{eq:A}).

\medskip
\STATE function $N$={\tt MPA}($x$,$t$)

\medskip
\STATE $B$={\tt TNBDLR}($x$,$t$);
\STATE $[Q,C]={\tt TNQR}(B)$;
%\STATE $\left(\begin{array}{l} d_1\\ d_2 \end{array}\right)=Q^Tb$;
%\STATE $z={\tt TNSolve}(R,d_1)$;
%\STATE $r=Q\left(\begin{array}{l} 0\\ d_2 \end{array}\right)$
\STATE $C_1=C(1:n+1,:)$;
\STATE $Q_1=Q(:,1:n+1)$;
%\STATE $d_1=Q_1'*b'$;
\STATE $IR_1={\tt TNInverseExpand}(C_1)$;
\STATE $MPA=IR_1*Q_1'$

\medskip
\STATE $N=MPA$;
\end{algorithmic}\label{alg:PseudoInverseA}
\end{algorithm}

Algorithm \ref{alg:PseudoInverseA} consists of computing the Moore-Penrose inverse $A^{\dagger}$ of the STP Lagrange-Vandermonde matrix $A$ in (\ref{eq:A}), starting from the matrix $\mathcal{BD}(A)$ computed to HRA by means of the algorithm ${\tt TNBDLR}$ described in Section \ref{sec:TNBDLR}. Let us observe that since $R$ is TP, which is consequence of the strict total positivity of $A$, its inverse $R^{-1}$ can be computed to HRA in $O(n^2)$ arithmetic operations by using the algorithm {\tt TNInverseExpand} in \cite{MM19} and whose implementation in {\sc Matlab} can be taken from \cite{KOEV}. As the cost of {\tt TNQR} is $O(l^2n)$, this is the total cost of Algorithm \ref{alg:PseudoInverseA}.

\begin{algorithm}[!h]
\caption{}
\begin{algorithmic}[1]
\REQUIRE Vectors $x\in \mathbb{R}^{n+1}$ and $t\in \mathbb{R}^{l+1}$ containing the nodes $\{x_j\}_{1\leq j\leq n+1}$ and $\{t_i\}_{1\leq i\leq l+1}$

satisfying (\ref{eq:NodesOrder}).
\ENSURE The Moore-Penrose inverse $L^{\dagger}$ of matrix $L$ in (\ref{eq:L}).

\medskip
\STATE function $N$={\tt MPL}($x$,$t$)

\medskip
\STATE $NA$={\tt MPA}($x$,$t$);

\medskip
\STATE {\bf for} $i=1:n+1$
\STATE \quad$den(i)={\tt prod}(x(i)-x(1:i-1))*{\tt prod}(x(i)-x(i+1:n+1))$;
\STATE \quad{\bf for} $j=1:l+1$
\STATE \quad\quad $MPL(i,j)=den(i)* NA(i,j)$;
\STATE \quad{\bf end}
\STATE {\bf end}

\medskip
\STATE $N=MPL$;
\end{algorithmic}\label{alg:PseudoInverseL}
\end{algorithm}

Algorithm \ref{alg:PseudoInverseL} corresponds to the calculation of $\overline D^{-1}A^{\dagger}$. The computational cost of the algorithm is leaded by the cost of function {\tt MPA} in line 2, which as we have explained is $O(l^2n)$.

Let us observe that the Moore-Penrose inverse $L^{\dagger}$ provides the solution of the least-squares problem we posed in Section \ref{sec:LS}. Explicitly, the solution of the least squares problem $Lc=b$, where $L$ is the matrix in (\ref{eq:L}), is the solution of the normal equations, $L^TL\overline c=L^Tb$. As the elements of $L$ satisfy (\ref{eq:NodesOrder}), matrix $L$ has full rank and the unique solution of the normal equations is given by $\overline c=L^{\dagger}b.$

\medskip

The {\em projection matrix} on the column space of matrix $L$ is $$H=LL^{\dagger}.$$ The next result proves that $H$ is also the projection matrix on the column space of matrix $A$ of (\ref{eq:A}).

\begin{proposition}
Given the nodes $\{x_j\}_{1\leq j\leq n+1}$ and $\{t_i\}_{1\leq i\leq l+1}$ such that the order given in  (\ref{eq:NodesOrder}) is satisfied, let $H$ be the projection matrix on the column space of matrix $L$ and $H_A$ be the projection matrix on the column space of matrix $A$ in (\ref{eq:L}). Then, it is satisfied that $H=H_A$.
\label{prop:ProjMatrix}
\end{proposition}

\begin{proof}
Since $L$ has full rank, by Proposition \ref{prop:MPInverse} it is satisfied that $L^{\dagger}=\overline D^{-1}A^{\dagger}$, and therefore $H=LL^{\dagger}=(A\overline D)(\overline D^{-1}A^{\dagger})=AA^{\dagger}=H_A.$
\end{proof}

Taking into account Proposition \ref{prop:ProjMatrix}, matrix $H$ can be computed from the $QR$ decomposition (\ref{eq:QR}) of $A$ as follows,
$$H=H_A=AA^{\dagger}=Q_1RR^{-1}Q_1^T=Q_1Q_1^T,$$
where $Q_1$ is the $(l+1)\times (n+1)$ matrix with the first $n+1$ columns of $Q$.

The next algorithm performs, taking advantage of the strict total positivity of matrix $A$, the computation of $H$ in $O(l^2n)$ arithmetic operations.

\begin{algorithm}[!h]
\caption{}
\begin{algorithmic}[1]
\REQUIRE Vectors $x\in \mathbb{R}^{n+1}$ and $t\in \mathbb{R}^{l+1}$ containing the nodes $\{x_j\}_{1\leq j\leq n+1}$ and $\{t_i\}_{1\leq i\leq l+1}$

satisfying (\ref{eq:NodesOrder}).
\ENSURE The projection matrix  $H$ of $L$ in (\ref{eq:L}).

\medskip
\STATE function $H$={\tt PL}($x$,$t$)

\medskip
\STATE $B$={\tt TNBDLR}($x$,$t$);
\STATE $[Q,C]={\tt TNQR}(B)$;
\STATE $Q_1=Q(:,1:n+1)$;
\STATE $PL=Q_1*Q_1'$;

\medskip
\STATE $PL$;
\end{algorithmic}\label{alg:ProjMatrix}
\end{algorithm}

Let us remark that, calling $\overline c$ to the solution of the least-squares problem $Lc=b$ we are dealing with, it is satisfied that the projection of $b$ onto the columns space of $L$ is $L\overline c=LL^{\dagger}b=Hb$. In consequence, matrix $H$ provides a way to evaluate the least-squares fitting polynomial at the nodes $\{t_i\}_{1\leq i \leq l+1}$.

\section{The general case}
\label{sec:GeneralCase}

Now we analyze the more general case in which the ordering of the nodes is not necessarily the specific one given in (\ref{eq:NodesOrder}). First we consider the least squares problem $Lc=b$, where $L$ is the matrix in (\ref{eq:L}).  In Proposition \ref{prop:SolLS} it was shown that such solution can be expressed as $\overline c=\overline D^{-1}\overline z$, where $\overline z$ is the least squares solution of the system $Az=b$, with $A$ being the Lagrange-Vandermonde matrix in (\ref{eq:A}). Let us observe that if (\ref{eq:NodesOrder}) is not satisfied, matrix $A$ is not STP and therefore Algorithm \ref{alg:LSLagMMV} cannot be used to solve this problem. In this situation, we can assume without loss of generality that $x_1<x_2<\ldots<x_{n+1}$ and that $t_{l+1}<\ldots<t_2<t_1$, and a change of variable
$$s=a_0+a_1t,~\text{ with }~a_0,a_1\in\mathbb{R},$$ can be done
  which takes all the $t_i$'s to a new interval to the right of $x_{n+1}$, so that (\ref{eq:NodesOrder}) holds and the least squares solution of $Az=b$ can be obtained from the least squares solution of the overdetermined linear system $My=b$, where $M$ is the STP Lagrange-Vandermonde matrix (\ref{eq:A}) at the new nodes $\{s_i\}_{1\leq i\leq l+1}$ as the next result shows.

\begin{proposition}
\label{prop:SolLSGeneral}
Given the nodes $\{x_j\}_{1\leq j\leq n+1}$ and the data $\{(t_i,b_i)\}_{1\leq i\leq l+1}$ such that $x_1<x_2<\ldots<x_{n+1}$ and $t_{l+1}<\ldots<t_2<t_1$, let $\overline y$ be the unique least-squares solution of the linear system $My=b$, where $M$ is the STP Lagrange-Vandermonde matrix (\ref{eq:A}) at the nodes $s_i=a_0+a_1t_i,~\text{ with }~a_0,a_1\in\mathbb{R},~1\leq i \leq l+1$. Then, the unique solution in the least squares sense of the linear system $Az=b$, where $A$ is the matrix in (\ref{eq:A}) at the original nodes $t_i$, is $\overline z=S^{-1}\overline y$, where $S$ is the matrix of change of basis from the basis $\mathcal{B}$ of (\ref{eq:basis}) to the basis

\begin{equation}
\label{eq:basisNewOrder}
\mathcal{B}_s=\Bigg\{ \prod_{\substack{k=1\\ k \neq 1}}^{n+1}(s-x_k), \prod_{\substack{k=1\\ k \neq 2}}^{n+1}(s-x_k), \ldots \prod_{\substack{k=1\\ k \neq n+1}}^{n+1}(s-x_k)\Bigg\},
\end{equation}
being $s=a_0+a_1t~(a_0,a_1\in\mathbb{R})$.
\end{proposition}

\begin{proof}
Solving in the least squares sense the system $Az=b$ consists of solving the normal equations $A^TA\overline z=A^Tb$. It can be seen that $A=MS$, and so solving the normal equations is equivalent to solve $S^TM^TMS\overline z=S^TM^Tb$. As $S$ is a matrix of change of basis, it is non singular and $M^TMS\overline z=M^Tb$. Calling $\overline y=S\overline z$ we have $M^TM\overline y=M^Tb$, whose solution $\overline y$ is the least squares solution of the overdetermined linear system $My=b$.
\end{proof}

This proposition together with Proposition \ref{prop:SolLS} provide the unique solution in the least squares sense of the overdetermined linear system $Lc=b$, which is given by $$\overline c=\overline D^{-1}S^{-1}\overline y.$$

We remark that although a change of variable is performed, the data vector $b$ does not change. 

Let us notice that for the computations we are interested in doing with the least squares fitting polynomial

\begin{equation}\label{eq:p}
p(t)=\overline c_1 l_1(t)+\overline c_2 l_2(t)+\cdots+\overline c_{n+1} l_{n+1}(t),
\end{equation}
where $\overline c=(\overline c_1 ~\overline c_2 \cdots \overline c_{n+1})^T$ and $l_i(t)$ are the Lagrange polynomials in (\ref{eq:basisLagrange}), the explicit computation of $\overline c$ and in consequence of $S$ and $S^{-1}$ are not required. Let us consider for example the problem of evaluating $p(t)$ at $t^*\in (t_{l+1},t_1)$. Writing $p(t)$ as a matrix product,

$$
p(t)=\overline c^T \left(
                     \begin{array}{c}
                       l_1(t) \\
                       l_2(t) \\
                       \vdots \\
                       l_{n+1}(t) \\
                     \end{array}
                   \right)=\overline c^T \overline D  \left(
                     \begin{array}{c}
                       l^*_1(t) \\
                       l^*_2(t) \\
                       \vdots \\
                       l^*_{n+1}(t) \\
                     \end{array}
                   \right),
$$
where $l^*_i(t)=\prod_{\substack{k=1\\ k \neq i}}^{n+1}(t-x_k)$, and taking into account that $\overline z=\overline D\overline c$, $A=MS$ and $\overline y=S\overline z$, it is obtained that
$$
p(t)=\overline z^T \left(
                     \begin{array}{c}
                       l^*_1(t) \\
                       l^*_2(t) \\
                       \vdots \\
                       l^*_{n+1}(t) \\
                     \end{array}
                   \right)=\overline z^TS^T\left(
                     \begin{array}{c}
                       l^*_1(s) \\
                       l^*_2(s) \\
                       \vdots \\
                       l^*_{n+1}(s) \\
                     \end{array}
                     \right)=\overline y^T\left(
                     \begin{array}{c}
                       l^*_1(s) \\
                       l^*_2(s) \\
                       \vdots \\
                       l^*_{n+1}(s) \\
                     \end{array}
                     \right),
$$
where $s=a_0+a_1t$. Thus, the evaluation of the least squares fitting polynomial whose coefficients $\overline c$ are the least squares solution  of the system $Lc=b$ at $t\in [t_{l+1},t_1]$ is equal to the evaluation of the least squares fitting polynomial whose coefficients $\overline y$ are the least squares solution of the system $My=b$ at $s=a_0+a_1t$.

\bigskip

As for the special case of evaluating the least squares fitting polynomial $p(t)$ in (\ref{eq:p}) at the points $\{t_i\}_{1\leq i\leq l+1}$, the projection matrix $H$ provides an alternative way to do it. By Proposition \ref{prop:ProjMatrix}, $H=H_A$ and, since $A=MS$ (where $S$ is nonsingular) proceeding analogously as in the proof of Proposition 3, we get that $H_A=H_M$, where $H_M$ is the projection matrix on the column space of matrix $M$. As $M$ is an STP matrix it follows that $H_M$, or what is the same $H$,  can be computed accurately by using Algorithm \ref{alg:ProjMatrix}. In consequence, the evaluation of $p(t)$ at  $\{t_i\}_{1\leq i\leq l+1}$ is performed accurately in this way
\begin{equation}
\left(
  \begin{array}{c}
    p(t_1) \\
    p(t_2) \\
    \vdots \\
    p(t_{l+1}) \\
  \end{array}
\right)=H b,
\label{eq:VectorProj}
\end{equation}
also without requiring the explicit computation of $S$ or $\overline c$.

Concerning the practical evaluation at points different from the points $t_i$ it must be observed that, as indicated in Chapter 5 of \cite{T13} (based on the works \cite{H04,WTG12}), the first form of the barycentric formula has better stability properties than the second form. In the context of our work, the adaptation of the first form algorithm to the case of a Lagrange-type basis (without denominators) reduces the complexity from $O(n^2)$ to $O(n)$.

\section{Numerical Experiments} %%%%%%%%%%%%%%%%%%%%%%%%%%%%%%%%%%%%%%%%%%%%%%%%%%%%%%%%%%%%%%%%%%%%%%%%%%%%%%%%%%%%%%%%%%%%%%%
\label{sec:NumExp}

Two numerical experiments illustrating the good performance of the algorithms developed in the previous sections are presented. In the first example, the nodes $\{x_j\}_{1\leq j \leq n+1}$ and $\{t_i\}_{1\leq i \leq l+1}$ have been chosen at random, with no other restriction than that of satisfying the ordering given in (\ref{eq:NodesOrder}). In addition, two data vectors are considered, the components of the first one being allowed to have different signs and all the components of the second one constrained to be positive. The second example corresponds to the general case, described in Section \ref{sec:GeneralCase}. The nodes $\{x_j\}_{1\leq j \leq n+1}$ and $\{t_i\}_{1\leq i \leq l+1}$ are Chebyshev points of the second kind which do not satisfy the ordering given in (\ref{eq:NodesOrder}), and the data vector $b$ is formed by the values of a given function at the  $t_i$. 

\begin{example}

Let $L$ be the collocation matrix corresponding to the Lagrange basis $\mathcal{B}_L$ of the polynomial space $\Pi_{20}(t)$ associated to the nodes

\bigskip

$$\begin{array}{ll}
\{x_j\}_{1\leq j \leq 21}~=&\hspace*{-0.1cm}\left\{-14,-\dfrac{129}{10},-\dfrac{116}{10},-11,-\dfrac{97}{10},-\dfrac{84}{10},-\dfrac{79}{10},-\dfrac{72}{10},-\dfrac{69}{10},-\dfrac{64}{10},-6,-5,-\dfrac{44}{10},-\dfrac{35}{10},\right.\\
&\left.~-3,-\dfrac{3}{2},-\dfrac{9}{10},-\dfrac{4}{10},0,\dfrac{1}{10},\dfrac{23}{100}\right\},
\end{array}$$

\noindent and

$$
\begin{array}{ll}
\{t_i\}_{1\leq i \leq 31}~=&\hspace*{-0.1cm}\left\{\dfrac{99}{10},\dfrac{96}{10},\dfrac{92}{10},9,\dfrac{87}{10},\dfrac{84}{10},\dfrac{81}{10},8,\dfrac{775}{100},\dfrac{75}{10},\dfrac{725}{100},7,\dfrac{68}{10},\dfrac{63}{10},6,\dfrac{59}{10},\dfrac{56}{10},\dfrac{52}{10},5,\dfrac{45}{10},\dfrac{41}{10},\dfrac{37}{10},\right.\\[5mm]
&\left.~\dfrac{32}{10},3,\dfrac{28}{10},\dfrac{23}{10},\dfrac{21}{10},\dfrac{16}{10},\dfrac{125}{100},1,\dfrac{8}{10}\right\}.
\end{array}$$

The condition number of $L$, computed in {\it Mathematica} by dividing the greatest by the smallest singular value of $L$, is $\kappa_2=4.1e+32$.

The two least squares problems associated to the two following data vectors will be solved,
$$\begin{array}{ll}
b^{^{(1)}}~=&\hspace*{-0.1cm}\left\{\dfrac{39}{100}, -\dfrac{17}{10}, -\dfrac{58}{10}, -4, 5, -\dfrac{57}{10}, \dfrac{63}{10}, -\dfrac{88}{100}, -\dfrac{39}{10}, \dfrac{69}{10},-7, \dfrac{44}{10}, 3, \dfrac{62}{10}, -\dfrac{57}{10}, -\dfrac{45}{10}, \dfrac{48}{10},-\dfrac{85}{100}, \right.\\[4mm]
&\left. ~ \dfrac{48}{10},\dfrac{24}{10},-4,\dfrac{28}{10}, \dfrac{27}{10}, \dfrac{46}{100}, -\dfrac{27}{10}, -\dfrac{12}{10}, -\dfrac{11}{10}, -\dfrac{15}{10}, \dfrac{12}{10},-\dfrac{84}{100}, -\dfrac{12}{100}\right\},\end{array}$$
$$\begin{array}{ll}
b^{^{(2)}}~=&\hspace*{-0.1cm}\left\{\dfrac{85}{10}, \dfrac{81}{10}, \dfrac{76}{10}, \dfrac{74}{10}, \dfrac{70}{10}, \dfrac{67}{10}, \dfrac{64}{10}, \dfrac{63}{10}, \dfrac{60}{10}, \dfrac{58}{10},\dfrac{56}{10}, \dfrac{53}{10}, \dfrac{51}{10}, \dfrac{47}{10}, \dfrac{45}{10}, \dfrac{44}{10}, \dfrac{41}{10}, \dfrac{38}{10}, \dfrac{37}{10}, \dfrac{33}{10},\dfrac{31}{10}, \dfrac{28}{10}, \right.\\[4mm]
&\left. ~ \dfrac{25}{10},\dfrac{23}{10}, \dfrac{22}{10},\dfrac{19}{10}, \dfrac{18}{10}, \dfrac{15}{10}, \dfrac{13}{10}, \dfrac{11}{10},\dfrac{97}{100}\right\}.\end{array}$$

\bigskip

\noindent Let $\overline{c}_i$, be the exact solution (in the least squares sense) of $Lc=b^{^{(i)}}$, computed with exact arithmetic in {\it Mathematica}, and let $c_i^*$ be an approximated solution of the problem, with $1\leq i \leq 2$. The relative error of $c_i^*$ is computed by using
\begin{equation}\frac{||c_i^*-\overline{c}_i||_2}{||\overline{c}_i||_2},\label{eq:RelErrLS}\end{equation}

\noindent and the results obtained with our approach (Algorithm \ref{alg:LSLagMMV}) and with the command $\backslash$ from {\sc Matlab} are shown in Table \ref{table:LSRelErrors}.

\begin{table}[!h]
\begin{center}
\begin{tabular}{|c|c|c|}
%\cline{1-3} & \multicolumn{2}{|c|}{Least squares solution} \\ \cline{1-3}
\hline
&Algorithm \ref{alg:LSLagMMV} & {\tt L}$\backslash${\tt b} \\ \hline
$b^{^{(1)}}$ & 3.8e-16 & 4.6e+00      \\ \hline
$b^{^{(2)}}$ & 6.7e-15 & 4.6e+00      \\ \hline
\end{tabular}
\end{center}
\caption{Relative errors obtained when computing the least squares solution of $Lc=b^{^{(i)}}$.} \label{table:LSRelErrors}
\end{table}

\bigskip

Next, we consider the Moore-Penrose inverse computation of matrix $L$. We compare the  results obtained with Algorithm \ref{alg:PseudoInverseL} of Section \ref{sec:MPInverseProjMatrix} and with the command {\tt pinv} of {\sc Matlab} by using the expression
\begin{equation}\frac{||L^{\dagger^*}-L^{\dagger}||_2}{||L^{\dagger}||_2}, \label{eq:RelErrorPS}\end{equation}

\noindent where $L^{\dagger^*}$ is the Moore-Penrose inverse of $L$ computed in {\it Mathematica} with exact arithmetic. The relative errors are given in the first two columns of Table \ref{table:PSHRelErrors}.

\begin{table}[!h]
\begin{center}
\begin{tabular}{|c|c||c|c|}
\cline{1-4} \multicolumn{2}{|c||}{Moore-Penrose inverse} & \multicolumn{2}{c|}{Projection matrix} \\ \cline{1-4}

Algorithm \ref{alg:PseudoInverseL} & {\tt pinv} & Algorithm \ref{alg:ProjMatrix} &  $L\cdot${\tt pinv}$(L)$ \\ \hline
1.9e-15 & 1.0e+00 & 1.8e-15 & 1.0e+00  \\ \hline
\end{tabular}
\end{center}
\caption{Relative errors obtained when computing the Moore-Penrose inverse and the projection matrix of $L$.} \label{table:PSHRelErrors}
\end{table}

Finally, the computation of the projection matrix $H$ of $L$ is carried out by our method (Algorithm \ref{alg:ProjMatrix}) and by multiplying $L$ and ${\tt{pinv}}(L)$ in {\sc Matlab}. The relative errors are calculated analogously as in (\ref{eq:RelErrorPS}), where the exact matrix $H$ is computed in {\it Mathematica} with infinite precision. The results are included in the last two columns of Table \ref{table:PSHRelErrors}.

Looking at the results shown in Table \ref{table:LSRelErrors} and Table \ref{table:PSHRelErrors} we observe that, although the 2-norm condition number of matrix $L$ is high, the algorithms developed in this work provide very accurate results, whereas generic algorithms of {\sc Matlab}, which do not take into account the specific structure of matrix $L$, give no accurate results at all.
\end{example}

\begin{example}

In this example $L$ is the collocation matrix corresponding to the Lagrange basis $\mathcal{B}_L$ of the polynomial space $\Pi_{11}(t)$ associated to the nonuniform sets of Chebyshev points of the second kind

$$\{x_j\}_{1\leq j \leq 11}=\left\{-\cos\left(\dfrac{(j-1)\pi}{10}\right)\right\}_{1\leq j \leq 11},$$
and
$$\{t_i\}_{1\leq i \leq 21}=\left\{\cos\left(\dfrac{(i-1)\pi}{20}\right)\right\}_{1\leq i \leq 21}.$$
The data corresponding to the $t_i$ are
$$\{b_i\}_{1\leq i \leq 21}=\{f(t_i)\}_{1\leq i \leq 21},$$
where
$$f(t)=\text{e}^t\sin(15t).$$
We observe that the Chebyshev points are not properly ordered, i.e., they do not fulfil (\ref{eq:NodesOrder}). Instead, the following relations are satisfied
$$-1=x_1<x_2<\cdots<x_{11}=1,$$
and
$$-1=t_{21}<t_{20}<\cdots<t_1=1.$$
In this situation the Lagrange-Vandermonde matrix $A$ of (\ref{eq:A}) is not STP, and a change of variable is performed so that the ordering given in (\ref{eq:NodesOrder}) is satisfied. In particular, we will perform the translation
$$s=\dfrac{11}{5}+t,$$
which takes the Chebyshev points $t_i$ (in decreasing ordering from 1 to $-1$) to the points $s_i$ (in decreasing ordering from $3.2$ to $1.2$), and so one has
$$\underset{-1}{\underbrace{x_1}}<x_2<\cdots<\underset{1}{\underbrace{x_{11}}}<\underset{1.2}{\underbrace{s_{21}}}<s_{20}<\cdots<\underset{3.2}{\underbrace{s_1}},$$
which is the ordering given in (\ref{eq:NodesOrder}). %That is, the datum corresponding to each point $s_i$ is $f(t_i)$, with $1\leq i\leq 21$.

With this change of variable the collocation matrix $M$ corresponding to nodes $\{x_j\}_{1\leq j \leq n+1}$ and $\{s_i\}_{1\leq i \leq l+1}$ is STP. As shown in Section \ref{sec:GeneralCase}, the evaluation of the least squares fitting polynomials whose coefficients are the least squares solutions of the problems $Lc=b$ and $My=b$ at any point $t\in[-1,1]$ and $s=\dfrac{11}{5}+t$, respectively, coincide.

\begin{table}[!h]
\begin{center}
\begin{tabular}{|c|c||c|c|}
\cline{1-4} \multicolumn{2}{|c||}{Solution of $My=b$} & \multicolumn{2}{c|}{Projection vector} \\ \cline{1-4}

Algorithm \ref{alg:LSLagTypeMMV} & $M\backslash b$ & $H\cdot b$&  $M\cdot${\tt pinv}$(M)\cdot b$ \\ \hline
1.2e-15 & 6.0e-03 & 1.1e-15 & 1.8e-01  \\ \hline
\end{tabular}
\end{center}
\caption{Relative errors of the computed least squares solution of the system $My=b$ and of the computed projection vector.} \label{table:RelErrProjM}
\end{table}

To calculate the least squares solution of the system $My=b$ we apply Algorithm \ref{alg:LSLagTypeMMV}. In the left part of Table \ref{table:RelErrProjM} the relative error of such computed solution compared to the least squares solution computed in {\it Mathematica} with 100 digits is shown. Also, the relative error obtained when computing the solution of the system $My=b$ with the command $\backslash$ of {\sc Matlab} is given. Both relative errors have been calculated in {\sc Matlab} with a formula analogous to (\ref{eq:RelErrLS}).

In addition, the projection vector of (\ref{eq:VectorProj}) is computed by multiplying the projection matrix $H$ given by Algorithm \ref{alg:ProjMatrix} by the data vector $b$. Such projection vector is also computed by the product $M\cdot {\tt pinv}(M)\cdot b$, and the corresponding relative errors are shown in the right part of  Table \ref{table:RelErrProjM}.

As we have remarked in Section \ref{sec:GeneralCase} the function $f(t)$ does not have to be re-evaluated at the new points $s_i$ (this would not be possible if the given data are experimental data, not obtained by evaluating a function), as the data do not vary with the change of variable.

The results in Table \ref{table:RelErrProjM} show that our approach gives much more accurate results than standard algorithms in {\sc Matlab}.

\end{example}

\section*{Acknowledgements} %%%%%%%%%%%%%%%%%%%%%%%%%%%%%%%%%%%%%%%%%%%%%%%%%%%%%%%%%%%%%%%%%%%%%%%%%%%%%%%%%%%%%%%%%%%%%%%%%%%%%%%%%%%%%%%%%%%

The authors are members of the Research Group {\sc asynacs} (Ref.CT-CE2019/683) of Universidad de Alcal\'a.

\bigskip
\bigskip
\bigskip
\bigskip
\bigskip

\bigskip
\bigskip

\end{document}